\documentstyle{amltd}
\begin{document}
\annalsline{151}{2000}
\received{January 8, 1999}
 \revised{Appendix November 30, 1999}
\startingpage{359}

\def\bye{\end{document}}
\font\smet=cmr7
\def\smallet{\hbox{\smet \'et}}

\def\ritem#1{\item[{\rm #1}]}

\def\k{\overline k}
\def\Z{\overline Z}
\def\C{\overline C}
\def\X{\overline X}
\def\U{\overline U}
\def\P{\overline P}
\def\V{\overline V}
\def\ra{\rightarrow}
\def\la{\leftarrow}
\def\Q{{\bf Q}}
\def\Spec{ { \rm Spec      }}

\input amssym.def
\input amssym.tex

 \title{Rational connectedness and Galois covers\\ of the projective line}
\shorttitle{Rational connectedness and Galois covers}  
\acknowledgements{}
 \author{Jean-Louis Colliot-Th\'el\`ene}
 \institutions{C.N.R.S., U.M.R. 8628, Universit\'e de Paris-Sud, Orsay, France
\\
{\eightpoint {\it E-mail address\/}:  colliot@math.u-psud.fr}}
 
 \vfill
Let $k$ be a $p$-adic field.
Some time ago, D. Harbater [9] proved that any finite group $G$
may be realized as a regular Galois group over the rational
function field in one variable $k(t)$, namely  there exists
a finite field extension $F/k(t)$, Galois with group $G$,
such that $F$ is
a regular extension of $k$ (i.e.\ 
$k$ is algebraically closed in $F$). Moreover, one may arrange  that a
given $k$-place of $k(t)$ be totally split in $F$. Harbater proved this
theorem for $k$ an arbitrary complete valued field.  Rather formal
arguments ([10, \S 4.5];  \S2 hereafter)
then imply that the theorem holds over any `large' field  $k$.
This in turn is a special case of a result of Pop [15], hence will be
referred to as the Harbater/Pop theorem.
We refer to
[10], [16], [6] for precise references to the literature  (work of
D\`ebes, Deschamps, Fried, Haran,
Harbater, Jarden, Liu,  Pop,  Serre, and
V\"olklein).

Most  proofs
(see [10], [19,  8.4.4,
p.~93] and  Liu's contribution to [16];  see however  [15])
first use direct arguments to establish the theorem
when $G$ is a cyclic group
(here the nature of the ground field is
irrelevant),  then proceed by patching, using either
formal or rigid geometry, together with
GAGA theorems.

In the present paper,  where I take the case of
algebraically closed  fields for granted, I show how a technique recently
developed
by  Koll\'ar [12] may be used to give a
quite different proof of the Harbater/Pop
theorem, when the `large' field $k$ has characteristic zero.
This proof  actually gives  more than the original result (see comment after
statement of Theorem 1).

Before I formally state the main result, let us recall what a `large' field is.
Let $k$ be a field and let $k((y))$ be the quotient field of the ring
$k[[y]]$ of formal power series in one variable.
Following F. Pop, we shall say that $k$ is `large'
if it satisfies one of the three equivalent properties  ([15,  Prop.~1.1]):
\begin{itemize}
\item[(i)] It is  existentially closed in $k((y))$: any  $k$-variety with a
$k((y))$-point has
a $k$-point.
 \item[(ii)] On a smooth integral $k$-variety with a $k$-point,
$k$-points are Zariski dense.

 \item[(iii)] On a smooth integral $k$-curve with a $k$-point,
$k$-points are Zariski dense.
\end{itemize}
(Such a field is clearly infinite. By going over to the completion
at a smooth $k$-point of a curve, one sees that (i) implies (iii).
That (iii) implies (ii) is easy (consider a regular system of
parameters). In characteristic zero, one may use resolution of
singularities to  show that (ii) implies (i).)

  Known examples of  `large' fields $k$  are
 fraction fields of a henselian discrete valuation ring,
 such as a $p$-adic field or a field of the shape
 $k=F((x))$ for $F$ some field.

Other well-known examples are real closed fields.
That these are `large' is a special instance of the following fact,
which seems to have escaped the attention of specialists:
any field $F$, all finite field extensions of  which
 are of degree a power of a fixed prime $p$,
is a `large' field. To see this,
one only needs to observe that on a regular, projective,
connected curve $C$
over a field $F$,
 given any nonempty open set $U$, any
 zero-cycle (divisor) $z$ on $C$  is rationally equivalent
 to a zero-cycle $z_1$ whose support is contained in $U$
 (a semi-local Dedekind ring is a principal ideal domain);
 the degree (over $F$) of $z$ and $z_1$ clearly coincide.
 Applying this to an $F$-point of $C$, one produces
 a zero-cycle $\sum_in_iP_i$ ($n_i \in {\bf Z}$, $P_i$ closed points)
 with support in $U$,
 such that the degree $\sum_i n_i [F(P_i):F]=1$.
 For $F$ as above, this forces one of the  degrees
 $[F(P_i):F]$ to be one.

 Other known examples are
 the fields of totally real algebraic numbers and
 of totally $p$-adic algebraic numbers (that these fields
 are `large' is a very special case of a theorem of Moret-Bailly
 [14, Thm.~1.3]).
 The property trivially holds for so-called pseudo algebraically
closed fields,
 such as infinite algebraic extensions of a finite field.

\nonumproclaim{Theorem 1} 
 Let $G$ be a finite group{\rm .} Let $k$ be a {\rm `}large{\rm '} field
of characteristic zero{\rm .}
Let ${\cal E}={\rm Spec}(K)$ be a $G$\/{\rm -}\/torsor
over ${\rm Spec}(k)${\rm .}
Then there exist an open set $U$
of the affine line  ${\bf A}^1_k$
containing a $k$\/{\rm -}\/point $O$ and
a $G$-torsor $V \rightarrow U$ such that the following two properties hold\/{\rm :}\/
 \begin{itemize}
\ritem{(i)} The fibre of $V \rightarrow U$
over $O$ is isomorphic to ${\cal E}$ {\rm (}\/as a $G$\/{\rm -}\/torsor
over ${\rm Spec}(k)${\rm );}
 
\ritem{(ii)} The smooth $k$\/{\rm -}\/curve $V$ is geometrically connected{\rm .}
\end{itemize}

\endproclaim

The ring $K$ is a finite separable extension of
$k$; it need not be a field.
In loose terms: given a Galois extension $K/k$ with group
$G$, one may realize $G$ as the Galois group of a `regular'
extension of $k(t)$, in such a way that over a suitable $k$-place of
$k(t)$, the extension specializes to $K/k$.

When  the $G$-torsor ${\cal E}/{\rm Spec}(k)$ is trivial,
i.e.\  ${\cal E}=\coprod_{g \in G} {\rm Spec}(k)$,
we recover the result of Harbater and Pop.
The question whether ${\cal E}$ may be chosen arbitrary
had been investigated for special groups
by several authors (see [6]).
For arbitrary groups, D\`ebes  proves
a weaker result ([6, Thm.~3.1])
when $k$ is `large', and he proves
the theorem  in the  case where $k$ is a pseudo
algebraically closed field ([6, Thm.~3.2]).

Using general results from  [EGA IV${}_3$],  we immediately
obtain a series of concrete corollaries. These will be detailed in Section 2.
 In the case of a split ${\cal E}/k$, most of them had already been obtained,
 with somewhat different proofs.

After the paper was submitted,  I was asked whether
in Theorem 1 one may impose arbitrary
$G$-torsors  as fibres of
$V \rightarrow U$  at more than one
$k$-point of $U \subset {\bf A}^1_k$.
The answer is in general in the negative, as shown in  the appendix.

Let us say a few words on the tools used in this article.
In a series of  papers which appeared in
1992, Koll\'ar, Miyaoka and Mori
 developed a technique
which enables them, under some assumptions, to smooth
a tree of rational curves into a single rational curve
([13, Thm.~(2.1)]; see also [11, Chap.~II. 7, pp.~154--158] and  [5, \S 4.2]).
That work was over an algebraically closed field.
In his recent paper [12],
Koll\'ar extends the technique over `large' fields (e.g. local fields).
Under certain assumptions,
he manages to deform a set of conjugate
${\bf  P}^1$'s into a single ${\bf   P}^1$ defined
over the ground field.
 From this he gets the finiteness  of the set of $R$-equivalence classes
 on $k$-points of a geometrically rationally connected
variety defined over a  local field $k$.
That the key lemma of [12] precisely holds
 for `large' fields provided the
incentive for the present paper.

The proof I give for Theorem 1
starts from the classical fact that a finite group $G$ is
a  Galois group over $k(t)$ when $k$ is algebraically closed
of characteristic zero.
It then uses a natural versal model for a
$G$-torsor, and applies the deformation result
of [12] to (a smooth compactification of)
the base space of this $G$-torsor.
The proof uses the existence of such a smooth compactification,
but it avoids any consideration of the divisor at infinity:
there is no discussion of inertia groups at all.

The idea of using a versal model
of a $G$-torsor, originally due to E. Noether,
has come up a number of times in the literature,
notably in work of E. Fischer, D. Saltman [17],  F.  A.  Bogomolov [1];
see  [20] and [21] for further references.

\demo{Acknowledgement} I am  much indebted to J\'anos Koll\'ar for having\break
shown me his work [12] while in progress.
I  thank Pierre D\`ebes, David Harbater and Laurent Moret-Bailly
for their interest in my paper.
Proposition A.3  was found during a stay at M.S.R.I.,
Berkeley, in September, 1999.

\section{Proof of Theorem 1}

In this section, we shall assume that the ground field
$k$ (which is of characteristic zero) is uncountable.
The proof in the
countable case will be given in
Section  2.

Let  $\k$ be an algebraic closure of $k$. Given
a $k$-scheme $Z$, let us write ${\overline Z}=Z\times_k\k$.

(1) Let $G$ be a finite group and ${\cal E}/{\rm Spec}(k)$ a $G$-torsor.
Let us fix an embedding of
$G$ into some  general linear group ${\rm GL}_n$.
Here $G$ is viewed as a constant (split) $k$-group scheme,
${\rm GL}_n$ is the linear group over $k$
and $i : G \rightarrow {\rm GL}_n$ is a homomorphism
of $k$-group schemes. Let $U= {\rm GL}_n/G $
be the affine $k$-variety of `left classes'.
This is the affine $k$-scheme whose ring
is the ring of invariants for $G$ acting on
the ring $k[{\rm GL}_n]$.
The projection map  ${\rm GL}_n \rightarrow U$ makes
${\rm GL}_n$ into a right $G$-torsor $V$ over $U$.
 The left action of  ${\rm GL}_n$ on itself
induces a left action of ${\rm GL}_n$ on $U= {\rm GL}_n/G $
and the projection $V \ra U$
is equivariant for these
(left) actions.

Let us recall basic facts from  noncommutative \'etale cohomology.
Given any smooth affine $k$-group scheme $H$, and any
commutative $k$-algebra $A$,
we denote by
$H^1_{\smallet}(A,H)$
 the pointed cohomology set which classifies
(\'etale) (right) $H\times_kA$-torsors
over ${\rm Spec}(A)$ (up to  nonunique isomorphism).
Such torsors will simply be called $H$-torsors over $A$.
For any such $A$, there is an ``exact
sequence"
$$  V(A) \rightarrow U(A) \rightarrow
H^1_{{\smallet}}(A,G)
\rightarrow H^1_{{\smallet}}(A,{\rm GL}_n).$$
Let us  detail this sequence.
The map $V(A) \rightarrow U(A) $ is the obvious one;
it respects the (left) action of ${\rm GL}_n(A)$ on both sets.
The right $G$-torsor $V \rightarrow U$ defines an element
 $\xi \in H^1_{{\smallet}}(U,G)$.
 To an element  $\rho \in U(A)={\rm Hom}_k({\rm Spec}(A), U)$,
the map $U(A)
\rightarrow H^1_{{\smallet}}(A,G)$
 associates the class
 $\rho^*(\xi) \in H^1_{{\smallet}}(A,G)$
 of  the pull-back $\rho^*(V \ra U)$, which is a $G$-torsor over $A$.
Two points $x,y \in U(A)$ have the same image in
$H^1_{{\smallet}}(A,G)$
if and only if there exists $\alpha \in {\rm GL}_n(A)$ such that $\alpha.x=y$.
By Grothendieck's version of Hilbert's Theorem 90,
the set $H^1_{{\smallet}}(A,{\rm GL}_n)$ classifies projective
modules of
rank $n$ over
$A$. Thus if $A$ is semi-local, or if $A$ is a Dedekind ring with trival class
group, then $H^1_{{\smallet}}(A,{\rm GL}_n)$ is reduced to one
element,
and for any right $G$-torsor ${\cal T}$ over $A$ there exists an element
$\rho \in
U(A)$ such that ${\cal T}$ and  $\rho^*(V\ra U)$ are
isomorphic  $G$-torsors over $A$. In particular, there exists
a $k$-point $P \in U(k)$ such that the fibre $V_P$ of $V$ above $P$
is a $G$-torsor isomorphic to the given ${\cal E}/k$. We shall fix
such a $k$-point~$P$.

\medbreak
(2)
By classical results  (see [19, Chap.~6]),
we know that  $G$ is a `regular' Galois
group over
$\k(t)$.  In other words there exist a nonempty open set $W$ of the affine
line ${\bf
A}^1_{\overline k}={\rm Spec}(\k[t])$ and a $G$-torsor over $W$ whose
underlying
variety is integral. Let $A$ be the semi-local ring
of $\k[t]$ at $t=0$ and $t=1$, and let $S={\rm Spec}(A)$.
Let us abuse notation and  call 0, respectively 1, the points
of $S$ defined by $t=0$, respectively $t=1$.
Changing coordinates and semi-localizing produces
a $G$-torsor ${\cal T}$ over $S$ such that ${\cal T}$ is an integral scheme.

By (1),  there exists a
nonconstant
${\overline k}$-morphism $\rho : S \rightarrow \U$
such that the pull-back of the $G$-torsor $\V \ra \U$
under $\rho$  is isomorphic to the $G$-torsor ${\cal T}/S$.
Given any $\alpha \in {\rm GL}_n(A)$, the $G$-torsor
$(\alpha.\rho)^*(\V \ra \U)$ is $G$-isomorphic to
the $G$-torsor ${\cal T}$. In particular, it is
an integral scheme.

\medbreak
(3)
The action of ${\rm GL}_n(\k) $ on $\U(\k)$ is transitive; hence
the obvious action of
${\rm GL}_n(\k) \times {\rm GL}_n(\k)$
 on $\U(\k) \times \U(\k)$ is also transitive.
 Reduction of $A$ modulo $t$ and modulo $t-1$ induces
a surjective homomorphism   ${\rm GL}_n(A) \ra {\rm GL}_n(\k) \times
{\rm GL}_n(\k)$.
 Thus given  two  points $M, N  \in
\U(\k)$,  there exists $\alpha \in {\rm GL}_n(A)$
such that $\alpha.\rho \in \U(A)$ sends the point $t=0$ to $M$ and
the point $t=1$ to $N$.

\demo{{R}emark} One should compare the present
general  position argument with  `Kuyk's lemma' (see [20, Lemma 4.5]).
\enddemo

(4) Since char($k$)=0, by Hironaka's theorem, there
exist smooth, projective, geometrically integral
$k$-varieties $X_1$ and $X$, with $V$ open in $X_1$
and $U$ open in $X$, together with a $k$-morphism
$p: X_1 \rightarrow X$ extending
the  map $V \rightarrow U$ and
 inducing a $k$-isomorphism $V \simeq
p^{-1}(U)$.

\medbreak
(5)
 According to a theorem of Koll\'ar, Miyaoka and Mori
([13]; [11,\break Thm.~II. 3.11, p.~118]),
to the point  ${\overline P} \in \U(\k) \subset \X(\k)$
one may associate  countably many proper
subvarieties $V_i$ ($i\in I$)
of the smooth projective variety $\X$ such that if
$f: {\bf   P}^1_{\k} \ra \X$ is a nonconstant morphism,
$f(0)={\overline P}$ and the image
of $f$ is not contained in the union of the $V_i$'s, then
$f$ is free over $0 \in {\bf   P}^1_{\k}$.  By definition
(see [11, II. 3.1, p.~113]), this means that the coherent cohomology group
$H^1({\bf   P}^1_{\k},f^*T_{\X}(-2))$ vanishes (here $T_{\X}$ denotes the
tangent bundle of ${\X}$),
which amounts to the hypothesis that
in Grothendieck's decomposition of the vector bundle $f^*T_{\X}$
over  ${\bf   P}^1_{\k}$ as a sum of line bundles ${\cal O}_{{\bf
P}^1}(n_j)$,
we have $n_j > 0$ for each $j$ (this is the ampleness property
for the vector bundle $f^*T_{\X}$ on ${\bf   P}^1_{\k} $ , see [11, II.3.8,
p.~116]).

Since $k$ is uncountable, there exists  a point $Q \in \U(\k)$, $Q \neq
{\overline P}$,
which does not lie on any of the $V_i$'s
 (proof: use a generically finite projection to
projective space and induct on dimension).
By (3),  there exists $\alpha \in {\rm GL}_n(A)$
such that $\alpha.\rho \in \U(A)$ sends the point $t=0$ to ${\overline P}$ and
the point $t=1$ to $Q$. Since $X/k$ is proper,
the morphism $\alpha.\rho: S \ra \U$ extends to a (nonconstant) morphism
$f: {\bf   P}^1_{\k} \ra \X$.
The image of $f$ contains $\P$ and is not contained
in the union of the $V_i$'s, since this image contains $Q$.
By the quoted theorem ([11, II.3.11]), we conclude:

\medbreak
(5.1) {\it The vector bundle $f^*T_{\X}$ on ${\bf   P}^1_{\k}$ is ample}.
\medbreak
On the other hand, we have:
\medbreak
(5.2)  {\it The underlying variety of the $G$-torsor $f^*(\V \ra \U)$ over
$f^{-1}(\U)$
is integral}. 
\medbreak
Indeed, this follows from the same statement for the restriction of this
$G$-torsor over $S={\rm Spec}(A) \subset f^{-1}(\U)$, which was
pointed out at the end of (2).

\medbreak
(6) We have now reached the  situation studied in [12].
Starting from $f: {\bf   P}^1_{\k} \ra \X$ such that $f(0)=\P$
and $f^*T_{\X}$ is ample, Koll\'ar ([12, 3.2], I change notation) produces,
over the ground field $k$,
a smooth integral $k$-curve $C$ with a $k$-point
$O$, a smooth geometrically integral $k$-surface
$Z$ proper over $C$, together with a
$k$-morphism $h: Z \ra X$,
with the following properties:
\medbreak
(6.a) The projection $Z \ra C$ admits a $k$-section $\sigma: C \ra Z$
which by $h$ is mapped  to $P \in X$.
\medbreak
(6.b) The geometric fibre $Z_{\overline O}$ of $Z \ra C$ at the point $O$
is a comb\break $D+\sum_{i \in I}C_i$ on $\Z$ (here $I$ is a nonempty
finite set, the $C_i$'s are
the teeth of the comb, see [11, II.7.7, p.~156]),
each component of which is
a nonsingular curve of genus zero;
the map ${\overline h}: \Z \ra \X$
sends $D$ to $\P$ and induces on $C_i$
a conjugate  of $f: {\bf   P}^1_{\k} \ra \X$.
\medbreak
(6.c) Over any closed point $M$  of $C$ different from $O$, the fibre $Z_M$
of $Z \ra C$ is $k(M)$-iso\-morphic to the projective line
${\bf   P}^1_{k(M)}$: the fibre  is a smooth,
geometrically irreducible, projective curve of
genus zero over the residue field $k(M)$, and it contains the
$k(M)$-rational point
$\sigma(M)$.
 \medbreak

(7) Since the map ${\overline h}: Z_{\overline O} \ra \X$ is not constant
(because its restriction to any $C_i$ is not constant),
the closed set $h^{-1}(P) \subset Z$ is a proper closed set.
Thus, after shrinking $C$, we may assume: for no $M \in C$ is $h$ constant
on the fibre $Z_M$ (note that on any fibre $Z_M$,
$h$ assumes the value $h(\sigma(M))=P\times_kk(M)$).

Let $\Omega \subset Z$ be the inverse image of
$U$ under $h$.  Note that $\Omega$ contains $\sigma(C)$,
hence the composite map $\Omega \subset Z \ra C$
is surjective.
Let $\Omega_1 \ra \Omega$ be the
inverse image of the $G$-torsor $V \ra U$
under $h: \Omega \ra U$.
Let $M$ be a closed point in $C$. We shall show: {\it
For all but finitely many
$M\in C${\rm ,} the total space of the induced $G$\/{\rm -}\/torsor
$\Omega_{1,M} \ra \Omega_M \subset Z_M \simeq {\bf   P}^1_{k(M)}$
is a smooth geometrically integral $k(M)$\/{\rm -}\/variety}.

To prove this, it is enough to prove the
corresponding statement over $\k$.
For the rest of the proof of (7), to simplify notation,
let us set $k=\k$. Points $M$
will be $\k$-rational points on $C$.
For $M \neq O$, the (nonempty) variety $\Omega_M$
is smooth and connected and the
variety $\Omega_{1,M}$ is a finite \'etale cover
of $\Omega_M$, hence is smooth.
To prove that a given  $\Omega_{1,M}, M \neq O,$ is integral, it is
thus enough to show that it is connected.

The inverse image in $\Omega_1$
of $D\cap \Omega$ is a disjoint union of copies
$D_g$ ($g \in G$)
of $D\cap \Omega$, each with multiplicity one; by (5.2) and (6.b), for a
given $i \in I$ the
inverse image in
$\Omega_1$ of each $C_i \cap \Omega $ is a (smooth) {\it connected} curve,
which
 meets {\it each} $D_g$ ($g \in G$), since $C_i$ meets $D$ (see (6.b)). Thus
$\Omega_{1,O}$, which is the inverse image of $D+\sum_{i\in I}C_i$,
is a  {\it reduced connected} divisor on $\Omega_1$.

 That $\Omega_{1,M}$  is connected for all but finitely many $M \in C$
now follows from the general lemma (where $X$ and $Y$ have nothing to
do with the previous $Y$ and $X$), to be applied to
$X=\Omega_1$ and $Y=\Omega$:

\nonumproclaim{Lemma}  Let $C$ be a smooth{\rm ,} connected curve over
 an algebraically closed field $k${\rm ,} and let $O \in C(k)${\rm .}
 Let $X${\rm ,} $Y${\rm ,} $C$ be smooth varieties over
 $k${\rm ,} equipped with faithfully flat $k$\/{\rm -}\/morphisms $X \ra Y$ and $Y \ra C${\rm .}
 Assume that  the generic fibre of $Y \ra C$
is smooth and geometrically integral{\rm .} Assume that $X \ra Y$
is finite and {\rm \'{\it e}}tale{\rm .} Assume moreover that
the inverse image of $O$ under the composite map
$X \ra Y \ra C$ is a connected divisor on $X$
and is not a multiple divisor{\rm .}  Then there exists a finite set
$S$ of points of $C$ such that for $M \in C,  M \notin S${\rm ,}
the inverse image $X_M$ of $M$ under
the composite map
$X \ra Y \ra C$ is a smooth connected variety{\rm .}
\endproclaim

\demo{Proof}  
Note first that $X$ is connected. Indeed if it was not connected,
the finite \'etale cover $X \ra Y$ would break up into
a disjoint union of finite \'etale (hence faithfully flat)
covers $X_i \ra Y$, and the fibre of
$X \ra Y \ra C$ over $O$ would not be connected.
Thus $X$ is connected; since it is smooth, it is integral.
Let $D$ be the normalization of
$C$ in the function field of $X$. This is a smooth integral curve,
and the map $D \ra C$ is flat and finite.
Since $X$ is normal, the map $X \ra C$ factors through $D$.
The finite (\'etale) map $X \ra Y$ factors through
the scheme $Y\times_CD$. The scheme $Y\times_CD$ is integral,
because $C$ is its own normalization in $Y$, since we have
assumed that  the generic fibre of
 $Y \ra C$ is geometrically integral.
 The finite map
of integral varieties $X \rightarrow Y\times_CD$ is dominant, hence surjective
as a morphism of schemes (it need not be flat). In particular,
it is surjective on $k$-points (recall $k=\k$).
The projection map $Y\times_CD \ra D$ is faithfully flat,
since it is obtained by base change from the faithfully flat map
$Y \ra C$. In particular, $ Y\times_CD \ra D$ is surjective on $k$-points.
We conclude that $X \ra D$ is surjective on $k$-points.
But then the scheme-theoretic inverse image of $O \in C$
under the map $D \ra C$ must consist of one reduced point,
since the inverse image of $O$ under the composite map
$X \ra D \ra C$ is a connected divisor which is not multiple.
Since  $D \ra C$ is finite and flat, this implies that $D \ra C$ is
an isomorphism.  Thus the function field of $C$ is algebraically
closed in the function field of $X$, hence the generic fibre
of $X \ra C$ is a smooth geometrically integral variety.
 By [EGA   IV${}_3$, (9.7.7)] 
 this implies the same statement for
all fibres of $X \ra C$ away
from a proper closed subset of $C$.  \enddemo

\phantom{hi}

(8) We finally make use of the hypothesis that the field $k$ is `large.'
Since the curve $C$ has a $k$-rational point, namely $O$,
this hypothesis implies that there exists a $k$-point $M$ on
$C$ away from the finitely many points excluded in (7), such that the map
${\bf   P}^1_k  \ra X$ induced by $h$ on the fibre $Z_M \simeq {\bf   P}^1_k$
does what we want:  the inverse image of
the $G$-torsor
$V \ra U$  under the map $h: h^{-1}(U) \cap {\bf   P}^1 \ra U$
is  a $G$-torsor over the  open set  $h^{-1}(U) \subset {\bf   P}^1_k$,
whose fibre at
 $\sigma(M) \in h^{-1}(U)(k) \subset {\bf   P}^1(k)$ is isomorphic to the
fibre of $V \ra U$
 at $P$, hence is isomorphic to  ${\cal E}$ (by the very choice of $P$, see
(1)), and whose  total space is a geometrically integral
$k$-variety (see (7)).

\section{Corollaries}

 \nonumproclaim{Theorem 2}
 Let  $O$ be a ${\bf Q}$\/{\rm -}\/point of the projective line
  ${\bf   P}^1_{\bf Q}${\rm .}
  Let $G$ be a
finite group and let
  ${\cal E}={\rm Spec}(K) \rightarrow
{\rm Spec}({\bf Q})$ be a $G$\/{\rm -}\/torsor{\rm .}
There exist a smooth{\rm ,} geometrically integral curve
$Y/{\bf Q}$ whose smooth compactification has  a ${\bf Q}$\/{\rm -}\/point{\rm ,}  an open
set
$U \subset {\bf   P}^1 \times_{\bf Q} Y$ containing $O\times_{\bf Q} Y${\rm ,}
and a $G$\/{\rm -}\/torsor $V \rightarrow U$ {\rm (}\/an {\rm \'{\it e}}tale Galois cover with group $G${\rm ),}
whose restriction to  $O \times_{\bf Q}  Y$ is the $G$\/{\rm -}\/torsor ${\cal
E}\times_{\bf Q}Y${\rm ,}
and
 such that the fibre of the composite
map $V \rightarrow U \rightarrow Y$ at any geometric point
of $Y$ is nonempty and connected {\rm (}\/hence integral{\rm ).}
\endproclaim

 \demo{Proof}  Let
$G \hookrightarrow {\rm GL}_{n,{\bf Q}}$ be an embedding.
The varieties $U,V,X,X_1$ which appear in the proof of Theorem 1
may all be defined  over ${\bf Q}$. We also have $P \in U({\bf Q}) \subset
X({\bf Q})$.

For any field $F$ with ${\bf Q} \subset F$, let us in this proof say that
an $F$-morphism $f: {\bf   P}^1_F \ra X_F$ is {\it good} if
$f(O)=P_F$ and  the inverse image of $V_F \ra U_F$
under $f$ (restricted to $f^{-1}(U_F)$) is a geometrically integral
$F$-variety.
Let $Z= {\rm Hom}_{\bf Q}({\bf   P}^1, X, O \mapsto P)$
(notation as in [11, II.1.4, p.~94]).
This is a countable union of ${\bf Q}$-varieties $Z_d$
($d$ for degree of the image of ${\bf   P}^1$,
in a fixed projective embedding of $X$).
An $F$-point of $Z$ will be called good if the corresponding
$F$-morphism $f: {\bf   P}^1_F \ra X_F$ is good.
Given arbitrary field extensions\break ${\bf Q} \subset E_1 \subset E_2$,
a point in  $Z(E_1)$ is good if and only if its image
in $Z(E_2)$ is good.

The field  ${\bf Q}((x))$ is uncountable.
By Theorem 1 over such a field,
as proved in Section 1, there exists a good ${\bf Q}((x))$-point on $Z$,
hence on $Z_d$ for some $d$.
Let $Y \subset Z_d$ be the scheme-theoretic closure
of the image of the corresponding morphism
${\rm Spec}({\bf Q}((x))) \ra Z_d$.
The ${\bf Q}$-variety $Y$ is geometrically integral.
We have the field embeddings ${\bf Q} \subset {\bf Q}(Y) \subset {\bf Q}((x))$.
Thus on the one hand the generic point of $Y$
is  a good ${\bf Q}(Y)$-point of $Z$;
on the other hand any\break ${\bf Q}$-compactification
of $Y$ has a ${\bf Q}$-point. Indeed, for any such compactification $Y_c$,
the map ${\rm Spec}({\bf Q}((x))) \ra Y$ extends to a ${\bf Q}$-morphism
${\rm Spec}({\bf Q}[[x]]) \ra Y_c$; the image of $x=0$ is a ${\bf Q}$-point of
$Y_c$.

Replacing  $Y$ by a nonempty open set,
one may ensure ([EGA IV${}_3$, (8.8.2)]) that
the corresponding good ${\bf Q}(Y)$-morphism
${\bf   P}^1_{{\bf Q}(Y)} \ra X_{{\bf
Q}(Y)}$
 extends to a\break $Y$-morphism
$\varphi: {\bf   P}^1 \times_{\bf Q} Y \ra X
\times_{\bf Q} Y$ which sends $O  \times_{\bf Q} Y$ to $P \times_{\bf Q} Y$.

Let $\Omega=\varphi^{-1}(U \times_{\bf Q} Y ) \subset {\bf   P}^1
\times_{\bf Q} Y$
and let $\Omega_1 \ra \Omega$ be the $G$-torsor which is the inverse
image of the $G$-torsor $V \times_{\bf Q} Y  \ra U \times_{\bf Q} Y $ under
$\varphi$.
Upon replacing $Y$ by a nonempty open set (this is actually not necessary),
the restriction of this $G$-torsor over $O
\times_{\bf Q} Y \subset \Omega$ is isomorphic to ${\cal E}\times_{\bf Q}Y$
(indeed, this is true over the generic point of $Y$).
We have the maps $\Omega_1 \ra \Omega \ra Y$.
The first map is finite \'etale of constant
rank, the second one is smooth and surjective.
Thus the composite map $\Omega_1 \ra Y$ is smooth.
 Since the
generic point of $Y$ corresponds to a good point of $Z$, the generic fibre
$\Omega_{1,{\bf Q}(Y)}$ is geometrically integral over ${\bf Q}(Y)$.
Upon replacing $Y$ by a nonempty open set ([EGA IV${}_3$, (9.7.7)(iv)]),
we therefore have that  all geometric fibres
of  the map $\Omega_1 \ra Y$ are smooth and geometrically integral.
In particular
 for any field $F$ with ${\bf Q} \subset F$ and any $F$-point of
$Y$, the morphism $\varphi_F: {\bf   P}^1_F \ra X_F$ induced by $\varphi$
is good.

On a smooth projective model $Y_c$ of $Y$ over ${\bf Q}$, there exists a
${\bf Q}$-point $R$.  By considering a regular system of
parameters at $R$ one produces a geometrically integral ${\bf Q}$-curve
$C
\subset Y_c$, smooth at $R$, and which meets $Y$. One now replaces $Y$ by
$Y \cap C$.
This completes the proof of Theorem 2. \enddemo

\phantom{hi}

\demo{Remarks and corollaries}

\smallbreak
(1) Note that  $Y$ in Theorem 2 need not have a
${\bf Q}$-point.  But for any field $k$ containing
${\bf Q}$ such that $Y(k) \neq \emptyset$,
$G$ is a `regular' Galois group over the rational field $k(t)$,
with the added information
that the fibre at the point $t=0$
is isomorphic to the torsor ${\cal E}\times_{\bf Q}k$.
This applies in particular to any `large' field of characteristic zero,
thus {\it completing the proof of Theorem {\rm 1} for fields which are countable}.
\smallbreak
(2) One should compare Theorem 2 with the contribution of  Deschamps in
[16], and the proof given here with that given in [7, 4.2].
\smallbreak
(3) One amusing corollary is that {\it for any finite group $G${\rm ,} there exists
a finite set of number fields $k_i$ such that the greatest common denominator of the
degrees $[k_i:{\bf Q}]$ is equal to one{\rm ,} and such that
$G$ is a {\rm `}regular{\rm '} Galois group over each $k_i(t)${\rm ,}
hence in particular a Galois group over
each $k_i$}.
The proof is simple:
on the smooth compactification $Y_c$ of the curve $Y$, there exists
a ${\bf Q}$-point, call it $M$. If we let $S \subset Y_c$ be the complement
of $Y$
in $Y_c$, there exists a zero-cycle $\sum_{i \in I}n_iP_i$ (here the $n_i$ are
integers, $P_i$ is a closed  point and $I$ is finite) on $Y_c$
which is rationally equivalent to
$M$, hence of degree one, and whose support is foreign to $S$, i.e.\  whose
support is
contained in $Y$.  Let $k_i$ be the residue field at the closed point $P_i$.
Then $\sum_{i \in I}n_i [k_i:{\bf Q}]=1$ and
$Y(k_i)\neq \emptyset$ for each $i$, hence the claim.

One could say that,
for any group $G$,
the inverse Galois group problem over ${\bf Q}$ acquires a positive answer
when passing from rational points to `zero-cycles of degree one.'

This could  have been noticed earlier.
For any prime $p$, let $K_p$ be the fixed field
of a pro-$p$-Sylow subgroup of the absolute Galois
group of ${\bf Q}$. As proved in the introduction
of this paper, $K_p$ is
 a `large' field.
By Theorem 1 (or, for that matter, the Harbater/Pop theorem),
$G$ is a regular Galois group  over
$K_p(t)$.
There  exists a finite subextension $L_p/{\bf Q}$ of $K_p/{\bf Q}$,
 such that $G$ is a regular Galois group  over $L_p(t)$.
 By Hilbert's irreducibility
theorem, $G$ is a Galois group over the number field $L_p$,
whose degree $[L_p:{\bf Q}]$ is prime to $p$.

\medbreak
(4) Starting from the statement of Theorem 2
and writing a model of the whole situation
over an open set of the ring of integers
(same references to [EGA IV${}_3$] as above),
one easily  deduces the following
result, which is a special case of a
theorem of Fried and V\"olklein:
{\it  For a given finite group $G${\rm ,} for almost all primes $p$
{\rm (``}\/almost all\/{\rm "} depending on $G${\rm ),} $G$ is a  {\rm `}\/regular\/{\rm '}
Galois group over ${\bf F}_p(t)$}  (see [10] and [7,
3.9] for references; in [7] a
model-theoretic argument  is given).
Simply note that if ${\cal Y}/{\bf Z}$ is a
smooth integral model of the smooth, geometrically integral curve
$Y/{\bf Q}$, then by classical estimates (Weil)
we have ${\cal Y}({\bf F}_p) \neq \emptyset$ for almost all primes $p$.
Here again, the present proof
enables us to get more: if we start off
with a given $G$-torsor ${\cal E}$ over a
nonempty open set of ${\rm Spec}({\bf Z})$,
 we may satisfy the additional requirement that
 for almost all primes $p$
the `regular' Galois extension over ${\bf F}_p(t)$
be unramified at  $t=0$,
the fibre being isomorphic to ${\cal E}\times_{\bf Z}{\bf F}_p$.
\enddemo

\def\diagram#1{\def\normalbaselines{\baselineskip=0pt\lineskip=5pt}
\matrix{#1}}

\def\hfl#1#2#3{\smash{\mathop{\hbox to#3{\rightarrowfill}}\limits
^{\scriptstyle#1}_{\scriptstyle#2}}}

\def\gfl#1#2#3{\smash{\mathop{\hbox to#3{\leftarrowfill}}\limits
^{\scriptstyle#1}_{\scriptstyle#2}}}

\def\vfl#1#2#3{\llap{$\scriptstyle #1$}
\left\downarrow\vbox to#3{}\right.\rlap{$\scriptstyle #2$}}

\def\ufl#1#2#3{\llap{$\scriptstyle #1$}
\left\uparrow\vbox to#3{}\right.\rlap{$\scriptstyle #2$}}

\def\pafl#1#2#3{\llap{$\scriptstyle #1$}
\left\Vert\vbox to#3{}\right.\rlap{$\scriptstyle #2$}}

\def\vl#1#2#3{\llap{$\scriptstyle #1$}
\left\vert\vbox to#3{}\right.\rlap{$\scriptstyle #2$}}

\bigbreak \centerline{\bf Appendix}

 \bigbreak

In this appendix, where for simplicity I assume
all  fields to be of characteristic zero,
I address the question:

\medbreak
 {\it Let $k$ be a field{\rm ,}   $G$ a finite group{\rm ,}
$n \geq 1$ an integer{\rm .} Let
 ${\cal E}_1, \cdots, {\cal E}_n$
be $G$\/{\rm -}\/torsors over $k${\rm .} Can one find
an open set $U \subset {\bf A}^1_k${\rm ,}
a $G$\/{\rm -}\/torsor $V \ra U$ and
$n$ points $P_1, \cdots, P_n \in U(k)$
such that for each $i${\rm ,}
the fibre $V_{P_i}$ is isomorphic to
${\cal E}_i$ as a $G$-torsor over} $k$?
\medbreak

Here are two cases where the answer is in the
affirmative:
\smallbreak
(i) $G$ is an abelian group, its $2$-primary subgroup
is of exponent $2^r$,  the cyclotomic field extension
$k(\mu_{2^r})/k$ is cyclic, and $n$ is arbitrary.
This is  a special case of [3, Thm.~7.9]
(various versions of this statement exist in
the literature; see  [17], [20]).
\smallbreak

(ii) $G$ is arbitrary, $k$ is `large'  and $n=1$: this is
Theorem 1 of the present paper (with the additional
piece of information that $V$ may be chosen
geometrically integral).
\smallbreak

In this appendix, I show by examples that for $n \geq 2$
and $k$ `large' the answer
to the above question is in general in the negative.

In the first part of the appendix, written in April 1999,
I consider the case left open in (i) above.  I  give an example
with  $G={\bf Z}/8$ and $k$ the\break 2-adic field ${\bf Q}_2$.
As may be expected, this example is closely related
to Wang's counterexample to Grunwald's theorem.

In the second part of the appendix, written in November 1999,
for an arbitrary  prime $p$, I give examples
with $G$ a  $p$-group and $k$ a suitable `large' field.
That part builds upon work of Saltman [18].

Background  and references for the first part of the
appendix (algebraic tori,
quasi-trivial and flasque tori,  groups of
multiplicative type, $R$-equivalence) will be found in
[2], [3], and [21]. For  $G$ a commutative algebraic group over  a  field $k$,
the \'etale cohomogy group $H^1_{{\smallet}}(k,G)$
may be identified with a Galois cohomology group, and will
be simply denoted $H^1(k,G)$.

\specialnumber{A.1} 
\proclaim{Proposition} Let $k$ be a field and
$A$ be a finite abelian group{\rm .}
One may embed the constant $k$\/{\rm -}\/group scheme $A$ into a commutative diagram
of exact sequences of $k$\/{\rm -}\/groups of multiplicative type\/{\rm :}
$$\diagram
{1  &  \ra   &  A         &  \ra  &   P_1      &  \ra  &    T           &
\ra   & 1  \cr
     &         & \vfl{}{}{2mm}  &         &  \vfl{}{}{2mm}  &        &
\vfl{}{=}{2mm}  &         &  \cr
1   &  \ra   &  F         & \ra    &  P_2       &  \ra   &  T         &
\ra  & 1 }
$$
where $T$ is a $k$\/{\rm -}\/torus{\rm ,} $F$ is a flasque  $k$\/{\rm -}\/torus
and $P_1$ and $P_2$ are quasi\/{\rm -}\/trivial  $k$-tori{\rm .}
\endproclaim

\demo{Proof}  By the well-known duality
$M \mapsto \hat{M}={\rm Hom}_{k-{\rm gr}}(M,{\bf
G}_{m,k})$ between
$k$-groups of multiplicative type and finitely generated Galois
modules over $k$, it is enough
to prove the dual result.
There exist exact sequences of finitely generated  Galois modules
$$ 0 \ra \hat{T} \ra \hat{P}_1 \ra \hat{A} \ra 0$$
and
$$ 0 \ra \hat{P} \ra \hat{F} \ra \hat{A} \ra 0$$
with $\hat{P}_1$ and $\hat{P}$ permutation modules, and $\hat{F}$ a flasque
module
(for the second sequence, see [3, (0.6.2)]). The pull-back
of the first sequence under the map $\hat{F} \ra \hat{A}$ is an exact sequence
$$ 0 \ra \hat{T} \ra \hat{P}_2 \ra \hat{F}  \ra 0$$
where the module $\hat{P}_2$  is an extension of the permutation module
$\hat{P}_1$
by the permutation module $\hat{P}$, hence is itself a permutation module.
Taking duals yields the proposition. \enddemo

For a quasi-trivial $k$-torus $P$, Hilbert's Theorem 90 implies $H^1(k,P)=0$.
Passing over to Galois cohomology in the  diagram of Proposition A.1, we get
the commutative diagram of exact sequences
$$\diagram{
P_1(k)     & \ra  & T(k)       & \ra  & H^1(k,A)  & \ra  &  0 \cr
 \vfl{}{}{2mm}  &        &  \vfl{}{=}{2mm}  &        &    \vfl{}{}{2mm}   &
&   \cr
P_2(k) & \ra & T(k) & \ra & H^1(k,F) & \ra  &  0  .\cr
}$$
From this diagram it immediately follows that the
map $H^1(k,A) \ra H^1(k,F)$ is onto.

Let us recall the following
basic fact from [2]: the map $T(k)  \ra  H^1(k,F)$
induces an isomorphism
$T(k)/R \simeq H^1(k,F)$. Here $R$ denotes $R$-equivalence
([2, \S 4]) on the set of $k$-points of the $k$-torus $T$.

\specialnumber{A.2}
\proclaim{Proposition}
With notation as above{\rm ,} assume that there exists
$\xi \neq 0 \in H^1(k,F)${\rm .} Let  $\eta \in H^1(k,A)$ denote
a lift of
$\xi$ under the surjective map $H^1(k,A) \ra H^1(k,F)${\rm .}
Then there do not exist an open set
 $U \subset {\bf A}^1_k$
 and an $A$\/{\rm -}\/torsor $X \ra U$ with the following properties\/{\rm :}
 there exist points $M, N \in U(k)$ such that
  the fibre of
 $X \ra U$ at $M$ is trivial while the  fibre of
 $X \ra U$ at $N$ has class $\eta \in H^1(k,A)${\rm .}
\endproclaim

\demo{Proof}  Let us assume there exist
such $U,M,N$. Since $P_1$ is a quasi-trivial $k$-torus,
for any $k$-scheme $V$ the \'etale cohomology group
$H^1_{{\smallet}}(V,P_1)$ is isomorphic to a sum of groups
${\rm Pic}(V\times_kK_i)$,  where the $K_i/k$
are finite separable field extensions of $k$.
For  $U \subset {\bf A}^1_k$, we thus have
$H^1_{{\smallet}}(U,P_1)=0$. Hence  the map
$ T(U) \ra H^1_{{\smallet}}(U,A) $
associated to the upper exact sequence in the diagram of
Proposition A.1 is onto. There thus exists
a $k$-morphism  $\varphi: U \ra T$
such that
$\varphi^*(P_1 \ra T)$ is isomorphic to the
$A$-torsor $X \ra U$.
The map $T(k) \ra H^1(k,A)$ sends
   $\varphi(M)$  to 0, and it
   sends $\varphi(N)$ to $\eta$.
Thus  the map $T(k) \ra H^1(k,F)$
sends  $\varphi(M)$ to 0, and it sends
 $\varphi(N)$  to $\xi \neq 0$.
Now since $U$ is an open set of
${\bf A}^1_k$, the
 points $\varphi(M) \in T(k)$
and $\varphi(N) \in T(k)$ are $R$-equivalent:
their images under the map $T(k) \ra H^1(k,F)$
should coincide. This contradiction establishes our contention. \enddemo

We still need to exhibit one case where the
hypotheses of Proposition A.2 are fulfilled.
Let $k$ be a field,
let $A={\bf Z}/8$ and let $T$ and $F$
be two $k$-tori as in Proposition A.1.
Suppose the  cyclotomic field extension $k(\mu_8)/k$
has degree 4. Its  Galois group is then
${\bf Z}/2 \times {\bf Z}/2$.
In that case, we have $H^1(k,{\hat F})={\bf Z}/2$ ([21, \S 7.4, p.~79]).
If $k$ is a $p$-adic field, then  the finite abelian groups
$H^1(k,S)$ and
$H^1(k,{\hat S})$ are dual (Tate-Nakayama).
Let $k$ be the 2-adic field ${\bf Q}_2$.
The field extension
 ${\bf Q}_2(\mu_8)/{\bf Q}_2$ has degree 4; we thus have
 $H^1({\bf Q}_2,F)\neq 0$.

 This completes the construction of the announced example,
 but
 one can  be more explicit.
 Let $k={\bf Q}_2$.
As a class  $\eta \neq 0 \in H^1(k,{\bf Z}/8) $,
let us take the class of the
degree 8 unramified field extension $E$ of $k={\bf Q}_2$.
Let us write the   commutative diagram in Proposition~A.1
over ${\bf Q}$. One may then write
the ensuing commutative diagram
over ${\bf Q}$ and over ${\bf Q}_2$,
in a compatible manner.
Let $M \in T(k)$ be any point with image
$\eta$ in  $H^1(k,{\bf Z}/8)$.
Suppose  the image of $\eta$  in $H^1(k,F)$ is trivial. Then
$M$ comes from a $k$-point of $P_2$. But then the
point $M$  lies in the closure of $T({\bf Q})$ in
$T({\bf Q}_2)$, since $P_2/{\bf Q}$ is a quasi-trivial torus,
hence  ${\bf Q}$-isomorphic to an open set of some affine space
over ${\bf Q}$. One can then find a
${\bf Q}$-point $N$ of $T$ such that the
fibre of $P_1 \rightarrow T$ at $N$ is
a Galois extension $F/{\bf Q}$ with group
${\bf Z}/8$ and such that $F \otimes_{\bf Q}{\bf Q}_2
\simeq E$ (as Galois extensions of ${\bf Q}_2$
with group ${\bf Z}/8$). But there is no such extension
(Wang's well-known counterexample to Grunwald's theorem,
see  [17] and [20]).
Thus the image of $\eta$  in $H^1(k,F)$ is nontrivial.

Let us now turn to other types of examples.

\specialnumber{A.3}
\proclaim{Proposition}  Let $p$ be a prime number{\rm .} There
exist a $p$\/{\rm -}\/group $G${\rm ,} a {\rm `}\/large\/{\rm '} field $k${\rm ,} and
$G$\/{\rm -}\/torsors ${\cal E}_1$ and ${\cal E}_2$ over $k$
with the following property\/{\rm :}
given any $G$\/{\rm -}\/torsor $f: V \ra U$ over an open set $U$
of  ${\bf A}^1_k${\rm ,} there do not exist
$k$\/{\rm -}\/points $P,Q \in U(k)$ such that the $G$\/{\rm -}\/torsor
$V_P$  is isomorphic to ${\cal E}_1$
and the $G$\/{\rm -}\/torsor
$V_Q$ is isomorphic to ${\cal E}_2${\rm .}
\endproclaim

\demo{Proof}    Saltman's work [18] (extended by Bogomolov [1],
see  [21, \S 7.6 and  \S 7.7]) produces finite $p$-groups $G$
together with faithful (finite dimensional) linear representations
$W$ of $G$ over the complex field ${\bf C}$, such that
the unramified Brauer group ${\rm Br}_{nr}(F)$ of $F={\bf C}(W)^G$
is a nontrivial ($p$-primary) group. Here by ${\bf C}(W)$ we denote the
fraction field of the symmetric algebra on $W$. The unramified Brauer
group of $F$ is the subgroup of the Brauer group ${\rm Br}(F)$
consisting of classes which are unramified with respect to
any (rank one) discrete valuation on $F$. As is well-known,
the group ${\rm Br}_{nr}({\bf C}(W)^G)$ does not depend on
the particular faithful (finite dimensional)
linear representation of $G$.

Let us fix one such
$p$-group $G$.
As in the beginning of Section 1, let us fix a homomorphic embedding
$G \ra {\rm GL}_n={\rm GL}_{n,{\bf C}}$. We may take for $W$ the vector space
of ${\bf C}$-points of $M_n$ (the ring scheme of $n$ by $n$ matrices over
${\bf C}$),
with
the action induced by left multiplication.
Let $U={\rm GL}_n/G$ and $V={\rm GL}_n \subset M_n$.
Projection $V \ra U$ makes $V$ into a $G$-torsor, whose properties are
described at the beginning of Section 1.

By Hironaka's theorem, there exists
 a smooth projective variety $X/{\bf C}$
containing
$U$ as a dense open set.
The function field ${\bf C}(X)$ of $X$ is $F$.
By results of Grothendieck, the natural map from
the \'etale Brauer group
${\rm Br}(X)=H^2_{{\smallet}}(X,{\bf G}_m)$
to ${\rm Br}(F)$  is one-to-one,  and it induces an isomorphism
 ${\rm Br}(X) \simeq {\rm Br}_{nr}(F)$ (see [4]).
 Let ${\cal A} \in {\rm Br}(X) \subset {\rm Br}(F)$
 be a nontrivial element. Let $X_F$ be the smooth, projective $F$-variety
 $X_F=X\times_{\bf C}F$. This contains the open set $U_F=U\times_{\bf C}F$.
On the one hand, the natural field embedding ${\bf C} \subset F$
  induces an inclusion
 $X({\bf C}) \subset X_F(F)$ of the
 set of ${\bf C}$-rational
 points of $X$ into the set of $F$-rational points of $X_F$, and similarly
 $U({\bf C}) \subset U_F(F)$. Let $P \in
 U_F(F)$ be an arbitrary  point in that subset.
On the other hand, the generic point ${\rm Spec}(F) \ra X$ of $X$
 gives rise (via the diagonal map) to an\break $F$-rational point $Q$ of $Y$.
Let ${\cal A}_F \in {\rm Br}(X_F)$ be the inverse image of ${\cal A}$ under
the projection map $X_F \ra X$. Let us evaluate ${\cal A}_F$
on the $F$-rational points $P$ and $Q$. We have ${\cal A}_F(P)=0 \in {\rm
Br}(F)$ because ${\cal A}_F(P)$ comes from ${\rm Br}({\bf C})$. We have
${\cal A}_F(Q) \neq 0 \in {\rm Br}(F)$ because ${\cal A}_F(Q)$ is none
other than the
image of ${\cal A} \in {\rm Br}(X)$ under the embedding ${\rm Br}(X)
\hookrightarrow
{\rm Br}(F)$. Let $k$ be a field, $F \subset k$, such
that the induced map ${\rm Br}(F) \ra {\rm Br}(k)$
is one-to-one. Changing the base field from $F$ to $k$,
we obtain  rational points which we still denote $P,Q$ in $X_k(k)$,
such that ${\cal A}_k(P)=0 $
and ${\cal A}_k(Q) \neq 0$ in ${\rm Br}(k)$. The points $P,Q$ both lie
in $U_k=U\times_{\bf C}k$.
Let ${\cal E}_1=V_P$, respectively
${\cal E}_2=V_Q$, be the $G$-torsors over $k$ defined
as the fibre of the $G$-torsor $V \ra U$ at $P$, respectively $Q$.
Suppose there exist a $G$-torsor $Z \ra Y$ over an open set
$Y \subset {\bf A}^1_k$ and two $k$-points $p, q \in Y(k)$
such that the fibre $Z_p$, respectively $Z_q$,  is a $G$-torsor over
$k$  isomorphic to ${\cal E}_1$, respectively ${\cal E}_2$.
By the general properties of the  $G$-torsor $V_k \ra U_k$ (see
beginning of \S 1)
and the fact that ${\rm Pic}(Y)=0$,
there exists a $k$-morphism $r: Y \ra U_k$ such that the inverse
image of the $G$-torsor  $V_k \ra U_k$ under $r$ is isomorphic to the
$G$-torsor $Z \ra Y$. Let $P_1=r(p) \in U(k)$ and $Q_1=r(q) \in U(k)$.
Then  $V_P$ and $V_{P_1}$ are  isomorphic as $G$-torsors
over $k$, and similarly $V_Q$ and $V_{Q_1}$.
The general properties
of the $G$-torsor $V \ra U$  then imply that there exist $g,h \in {\rm GL}_n(k)$
such that $gP_1=P$ and $hQ_1=Q$.  Since ${\rm GL}_n$ is an open set of an
affine space over $k$, this implies that the $k$-points
$P_1$ and $P$ of $U_k(k) \subset X_k(k)$ are $R$-equivalent. Similarly,
$Q_1$ and $Q$ are $R$-equivalent. Clearly, $P_1$ and $Q_1$
are $R$-equivalent. Thus $P$ and $Q$ are $R$-equivalent on
the projective $k$-variety $X_k$. By Prop.~16 of [2] (p.~213)
this implies  ${\cal A}_k(P)={\cal A}_k(Q)$.
But then we cannot have ${\cal A}_k(P)=0$ and ${\cal A}_k(Q) \neq 0$.

To complete the proof of Proposition A.3, it remains to
notice that the  field
$k=F((t))$ of formal power series in one variable
is a `large' overfield  of $F$ for which the  map
${\rm Br}(F) \ra {\rm
Br}(k)$ is one-to-one. \enddemo

Whether  examples as in Proposition A.3  may be exhibited
over a $p$-adic field  remains to be seen.

\references

[1]   
\name{F. A. Bogomolov}, The Brauer group of quotient
spaces of linear representations, {\it Izv.\ Akad.\ Nauk SSSR Ser.\ Mat\/}.\ {\bf 51}
(1987) 485--516;
Engl.\ transl.\ {\it Math.\ USSR Izv\/}.\ {\bf 30} (1988) 455--485.
 
[2]  
\name{J.-L. Colliot-Th\'el\`ene} et \name{J.-J. Sansuc},
La $R$-\'equivalence sur les tores, {\it Ann.\ Sci.\  \'Ecole Norm. Sup\/}., 4\`eme
s\'erie  {\bf 10}
(1977), 175--229.

[3]  
\bibline,
Principal homogeneous spaces
under flasque tori: applications,
 {\it J. of Algebra\/} {\bf 106} (1987), 148--205.

 [4]  
\bibline, The rationality problem
 for fields of invariants under linear algebraic groups
 (with special regards to the Brauer group), Lecture notes,
 IX ELAM, Santiago de Chile, 1988.

[5]  \name{O. Debarre}, Vari\'et\'es de Fano, in {\it S\'eminaire Bourbaki},
Vol.~1996-97, Exp.~827,
Ast\'erisque {\bf 245},
Soci\'et\'e Math.\  de France (1997), 197--221.

[6]  
\name{P. D\`ebes}, Galois covers with prescribed fibers: the Beckmann-Black
problem,
{\it Ann. Scuola Norm.\ Sup.\ Pisa Cl.\ Sci\/}.\   {\bf XXVIII} (1999) 273--286.

[7]  
\name{P. D\`ebes} and \name{B. Deschamps}, The regular inverse Galois
problem over large fields,
in {\it Geometric Galois Actions} (L. Schneps and P. Lochak, ed.),
London Math.\ Society Lecture Notes Series {\bf 243}  (1997), 119--138, 
Cambridge University Press, Cambridge.

[8]  
[EGA IV${}_3$]
\name{A. Grothendieck},
{\it  \'El\'ements de G{\rm \'{\it e}}om{\rm \'{\it e}}trie Alg{\rm \'{\it e}}brique},
r\'edig\'es avec la colla\-boration de J.~Dieudonn\'e.
IV.  \'Etude locale des
sch\'emas et des morphismes de sch\'emas (Troisi\`eme Partie),
Inst.\ Hautes \'Etudes Sci., Publ.\ Math.\ No.~{\bf 28}, 1966.

[9]  
\name{D. Harbater}, Galois coverings of the arithmetic line,  in
{\it Number Theory Seminar},  (New York 1984/85),  Lecture Notes in Math.\
{\bf 1240} (1987), 165--195, Springer-Verlag, New York. 

[10]  
\bibline, Fundamental groups of curves in characteristic
$p$, in {\it Proc.\  1994 Internat.\break  Congress of Mathematicians}
(Z\"urich) {\bf 1}, 656--666, Birkh\"auser Verlag, Basel,
 1995.  

[11]  \name{J. Koll\'ar}, {\it  Rational Curves on Algebraic Varieties},
Ergebnisse der Mathematik und ihrer Grenzgebiete, 3.\ Folge,
Band {\bf 32}, Springer-Verlag, New York,  1996.

[12]  
\bibline, Rationally connected varieties over local fields,
{\it Ann.\  of Math\/}.\  {\bf 150} (1999),  357--367.

[13] 
 \name{J. Koll\'ar, Y. Miyaoka},  and \name{S. Mori},  Rationally connected varieties,
{\it J. Algebraic Geom\/}.~{\bf 1}
(1992), 429--448.

[14] 
 \name{L. Moret-Bailly}, Groupes de Picard et probl\`emes de Skolem II,
{\it Ann.\ Sci.\  \'Ecole Norm. Sup\/}., 4\`eme s\'erie {\bf 22} (1989), 181--194.

[15] 
 \name{F. Pop}, Embedding problems over large fields, {\it Ann.\  of Math\/}.\ 
{\bf 144}
(1996), 1--34.

[16]  
{\it Recent developments in the inverse Galois problem},
Proc.\  of a 1993
Seattle research  conference,  
{\it Contemp. Math\/}. {\bf 186} (1995), 217--238.

[17]  
\name{D. J. Saltman}, Generic Galois extensions and problems in field theory,
{\it Adv.\ in Math}.\ {\bf 43} (1982),  250--283.

[18]  
\bibline,  Noether's problem over an algebraically
closed field, {\it Invent.\ Math}.\ {\bf 77} (1984), 71--84.

[19] 
 \name{J-P. Serre},  {\it Topics in Galois Theory}, Notes written by H.
Darmon,  Research Notes in Math., Jones and Bartlett
Publishers, Boston, MA, 1992.

[20]  
\name{R. G. Swan},  Noether's problem in Galois theory,
in {\it Emmy Noether in Bryn Mawr}, Proc.\ Symp.\
(Bryn Mawr, 1982) (Bh.\ Srinivasan and J. Sally, ed.),  Springer-Verlag, New York  (1983),
21--40.

 [21]  
\name{V. E. Voskresenski\v{\i}}, {\it Algebraic Groups and their Birational
Transformations\/}, Transl.\  of Math.\ Monographs {\bf 179},
A.M.S., Providence, RI, 1998.

\endreferences

\bye